\documentclass{amsart}

\begin{document}

\title{SYMMETRIES AND HIDDEN SYMMETRIES FOR FIELDS OUTSIDE BLACK HOLES\footnote{To appear in ``XVIIth International Congress on Mathematical Physics''.}}

\author{P. BLUE}

\address{
School of Mathematics, Maxwell Institute\\ University of Edinburgh,
Edinburgh, UK\\
Email: p.blue@ed.ac.uk\\
www.maths.ed.ac.uk/~pblue}


\maketitle

\begin{abstract}
This note surveys how energy generation and strengthening has been used to prove Morawetz estimates for various field equations in Minkowski space, the exterior of the Schwarzschild spacetime, and the exterior of the Kerr spacetime. It briefly outlines an approach to proving a decay estimate for the Maxwell equation outside a Kerr black hole. 
\end{abstract}


\section{Introduction}
This note crudely outlines a still-tentative approach to proving decay estimates for the Maxwell equations outside a slowly rotating Kerr black hole. This is part of a larger programme of research, undertaken by several research groups, to prove that the Kerr black holes are asymptotically stable against small perturbations. While it is physically unimaginable that a small ripple of gravitational radiation could cause a black hole to degenerate into a naked singularity or otherwise destroy the structure of a distant, asymptotically flat region, this is a challenging mathematical problem. The model for this approach is the original proof of the Cauchy stability of Minkowski space ($\mathbb{R}^{1+3}$ with $g=-\mathrm{d}t^2+\mathrm{d}x^2+\mathrm{d}y^2+\mathrm{d}z^2$)\cite{ChristodoulouKlainerman:MinkowskiStability}. That proof built on earlier work that used energies generated and strengthened by the (conformal) Killing vectors to study the wave equation, the Maxwell equation, and the linearised Einstein equation. 

\section{Energy estimates}
Let $M$ be a globally hyperbolic manifold, which is foliated by Cauchy hypersurfaces $\Sigma_t$. Assume further that the foliating hypersurfaces $\Sigma_t$ are the level sets of a function $t$, which we will refer to as the time. The contents of this section are well-known, and a relevant presentation can be found in \cite{ChristodoulouKlainerman:LinearFields}. Here, we emphasise the particle-wave analogies. By this, we mean the similarities in statements about null geodesics (representing massless particles) and about PDEs with a well-behaved energy-momentum tensor, such as the wave equation. 

Given a vector field $X$, the standard definition of the energy of a null geodesic, $\gamma$, with respect to $X$ and evaluated at time $t$ as 
\begin{equation*}
e_X[\gamma](t) = -\dot{\gamma}_\alpha X^\alpha
\end{equation*}
evaluated at the unique point where $\gamma$ intersects $\Sigma_t$. 

Classical field theories often have an energy-momentum tensor $\mathrm{T}_{\alpha\beta}$ that is symmetric ($\mathrm{T}_{\alpha\beta}=\mathrm{T}_{(\alpha\beta)}$), is divergence free ($\nabla^\alpha\mathrm{T}_{\alpha\beta}=0$), and satisfies the dominant energy condition (for all time-like, future-oriented vector fields $V$, $W$: $\mathrm{T}_{\alpha\beta}V^\alpha W^\beta\geq 0$). To be a little more careful, we should say a classical field theory gives a partial differential equation governing a (possibly indexed) field $\varphi$, that for each sufficiently smooth field $\varphi$ defined on $M$, there is a corresponding $\mathrm{T}[\varphi]_{\alpha\beta}$ which is a tensor field on $M$, and that $\mathrm{T}[\varphi](t)_{\alpha\beta}$ satisfies the above properties when $\varphi$ satisfies the relevant PDE. These properties will be assumed for the rest of the paper, although not every field theory satisfies these properties. The Einstein equation does not have a energy-momentum tensor, but there is a well known quadruple-indexed Bel-Robinson tensor with many of the same properties. 

Given an energy-momentum tensor as above, the standard definitions for the $4$-momentum $P_X[\varphi]_\alpha$ (of $\varphi$ with respect to a vector field $X$) and for the energy $E_X[\varphi](t)$ (of $\varphi$, generated by $X$, and evaluated on $\Sigma_t$) are
\begin{eqnarray*}
P_X[\varphi]_\alpha&=&\mathrm{T}_{\alpha\beta}X^\beta ,\\
E_X[\varphi](t)&=& \int_{\Sigma_t} P_X[\varphi]_\alpha \mathrm{d}\nu^\alpha ,
\end{eqnarray*}
where integration with respect to $\mathrm{d}\nu^\alpha$ denotes the standard flux integral through $\Sigma_t$ with respect to the normal vector $\nu$, which is assumed to be future-directed. 

There are two energy-generation properties: 
\begin{itemize}
\item{EG1:} If $T$ is a time-like and future-oriented vector field, then the energy of a null geodesics or a field is positive: 
\begin{eqnarray*}
e_T[\gamma]\geq0, &{}& E_T[\varphi]\geq0 .
\end{eqnarray*}
For null geodesics, the positivity results from taking the inner product of two causal, future-oriented vectors (and the minus sign included in the definition of the energy). For a field, the energy is the integral of the energy-momentum tensor evaluated on the generating vector $T$ and the hypersurface normal $\nu$, both of which are time like and future oriented. Thus, the integrand in the energy is positive by the dominant energy condition. 
\item{EG2:} Let $\gamma$ be parameterised by $t$. Let $\Omega[t_1,t_2]$ denote the region between $\Sigma_{t_1}$ and $\Sigma_{t_2}$, and let $\mathrm{d}\mu$ the spacetime volume element. The change in the energy between one time and the next is
\begin{align*}
e_X[\gamma](t_2)-e_X[\gamma](t_2)&=\int_{t_1}^{t_2} \dot\gamma_\alpha\dot\gamma_\beta 2\nabla^{(\alpha}X^{\beta)} \mathrm{d}t\\
E_X[\varphi](t_2)-E_X[\varphi](t_2)&=\int_{t_1}^{t_2} \mathrm{T}_{\alpha\beta} 2\nabla^{(\alpha}X^{\beta)} \mathrm{d}\mu.  
\end{align*}
The first property follows from integrating the derivative $\dot\gamma^\alpha\nabla_\alpha(\dot\gamma_\beta X^\beta)$ and applying the geodesic equation. The second property follows from the divergence theorem and the divergence-free property of $\mathrm{T}$. 
\end{itemize}
Noether's theorem can be seen as an application of the second energy-generation property: If $X$ is Killing, then $\nabla^{(\alpha}X^{\beta)}=0$ and the energies are conserved. Thus, if there is a smooth family of symmetries of the spacetime, there is a conserved energy. This is not the only application of the energy generation property. 

Conversely, symmetries are also relevant in other ways in the study of partial differential equations. For many geometric problems, the Lie derivative along a Killing vector gives a differential operator that commutes through the differential equations. In some cases, there might be a differential operator, $S$, which does not necessarily commute through the differential equation, but is a symmetry in the weaker sense that if $\varphi$ is a solution, then $S\varphi$ is also a solution. The existence of commuting operators is also related to the possibility of applying separation of variables in the analysis of the equation. 

Useful examples to consider are the linear $x$ momentum and the total angular momentum squared in Minkowski space. These are respectively quantities associated with $\partial_x$ and $\sum_{i=1}^3 \Theta_i^\alpha\Theta_i^\beta$ where $\Theta_i=\epsilon_{ij}{}^k x^j\partial_k$ are the vector fields generating rotations about the $x^i$ axis (with $\epsilon_{ijk}$ being the totally Levi-Civita symbol on constant $t$ hypersurfaces). Each of these tensors gives rise to a conserved quantity for null geodesics and to a differential operator that commutes with the d'Alembertian operator $\nabla^\alpha\nabla_\alpha$ and which has eigenfunctions that can be used in separation of variables (i.e. $e^{ikx}$ and the spherical harmonics $Y_{l,m}(\theta,\phi)$). In going to the Kerr case, it will be of particular importance to consider the case of the total angular momentum, for which the conserved quantity for null geodesics ($(\dot\gamma^\theta)^2+\sin^{-2}\theta(\dot\gamma^\phi)^2$) does not arise as the energy generated by any vector field and for which the commuting operator $\sin^{-1}\theta\partial_\theta\sin\theta\partial_\theta+\sin^{-2}\theta\partial_\phi^2$) does not arise as the Lie derivative along any vector field. We will discuss these properties in more detail in Sec. \ref{S:Kerr}.

For the discussion in this note, there is one key energy-strengthening property: 
\begin{itemize}
\item{ES:} If $S$ is a symmetry in the sense above, $X$ is a vector field, and $\varphi$ is the solution of some differential equation, then in addition to $E_X[\varphi]$, we can consider $E_X[S\varphi]$, which will enjoy properties EG1-EG2. If we are given a graded family of symmetries $\mathcal{S}=\cup_{n=0}^\infty\mathcal{S}_n$ where each symmetry in a given grade $S\in\mathbb{S}_n$ is a differential operator of order $n$, we can define the order $n$ energy with respect to $\mathbb{S}$ by
\begin{align*}
E_{X,n}[\varphi](t)=\sum_{i=0}^n\sum_{S\in\mathbb{S}_i}E_X[S\varphi](t) .
\end{align*}
\end{itemize}

\section{Some simple examples: geodesics and the wave equation in the Minkowski and Schwarzschild spacetimes}
The energy-momentum tensor for the wave equation $\nabla^\alpha\nabla_\alpha u=0$ is $\mathrm{T}[u]_{\alpha\beta}=\partial_\alpha u\partial_\beta u-(1/2)g_{\alpha\beta}(\nabla^\gamma u\nabla_\gamma u)$. In Minkowski space, there are the well-known energy and Morawetz estimates
\begin{eqnarray*}
E_T[u](0)=E_T[u](t)&=&\int_{\{t\}\times\mathbb{R}^3} \left(|\partial_t u|^2 +|\partial_r u|^2 +\frac{|\nabla_{S^2}u|^2}{r^2}\right)\mathrm{d}^3x  , \\
E_T[u](t_2)+E_T[u](t_1)&\geq& C\int_{[t_1,t_2]\times\mathbb{R}^3} \left(\frac{|\partial_t u|^2+|\partial_r u|^2}{r^2}+\frac{|\nabla_{S^2}u|^2}{r^3}+\frac{|u|^2}{r^4}\right) \mathrm{d}^3x\mathrm{d}t ,
\end{eqnarray*}
where $(t,r,\theta,\phi)\in\mathbb{R}\times(0,\infty)\times S^2$ are the standard spherical coordinates, $T=\partial_t$, and $C$ is some constant. The Morawetz estimate is a (quite weak) decay estimate, in the sense that $u$ and its derivatives must, on average, tend to zero as $t\rightarrow\infty$ in fixed regions of $r$ for the $t$ integrability property to hold. 

The positivity and conservation for the $T$ generated energy follow, respectively, from $\partial_t$ being time like (and future oriented) and being Killing. For a null geodesic, the same properties ensure that $e_T[\gamma](t)$ is positive and conserved. With $R=\partial_r$, we find
\begin{align*}
e_T&\geq e_R &&&
\dot\gamma_\alpha\dot\gamma_\beta\nabla^\alpha R^\beta&=\frac{2}{r^3}((\dot\gamma^\theta)^2+\sin^{-2}\theta(\dot\gamma^\phi)^2) .
\end{align*}
The first statement follows since the energy is linear in the generating vector field and since $T-R$ is causal and future directed. The second statement follows from direct computation. These two statements and the second energy generation property should give credibility to the Morawetz estimate with the time and radial derivative terms on the right removed. By replacing $R$ by $A=f(r)\partial_r$ (with $f(r)$ an appropriate choice of bounded, positive weight), it is possible to gain additional control over the remaining terms on the right of the Morawetz estimate. However, care must be taken in treating the boundary term at $r=0$ in any spherical coordinate calculation. In essence, the $A$ energy measures the radial momentum, which, along any null geodesic is always increasing. Thus, the terms appearing in the derivative of the $e_A$ energy are positive. 

In Minkowski space, translations parallel to the coordinate axes generate symmetries, so the partial derivative in these directions commutes through the d'Alembertian. These can be iterated to generate higher-order differential operators which are symmetries, and all combinations of these of a given order can be gathered into a set of symmetries of order $n\geq1$, $\mathbb{S}_n=\{\partial^{i_1}_x\partial^{i_2}_y\partial^{i_3}_z |i_1+i_2+i_3=n\}$. On each hypersurface of constant $t$, the (spatial) homogeneous Sobolev norm of order $n$ can be estimated by
\begin{align*}
\| u|_{\Sigma_t}\|_{\dot{{H}}^n}^2
\leq&\|u|_{\Sigma_t}\|_{\dot{{H}}^n}^2+\|(\partial_tu)|_{\Sigma_t}\|_{\dot{{H}}^{n-1}}^2 \\
&=\sum_{S\in\mathbb{S}_{n-1}}E_T[S u](t)=\sum_{i_1+i_2+i_3=n-1} E[\partial^{i_1}_x\partial^{i_2}_y\partial^{i_3}_z u](t) .
\end{align*}
(We must use the $n-1$ energy since the energy-momentum tensor has the derivative of $u$, instead of $u$ itself.) Control of such norms, and more complicated Klainerman-Sobolev norms, is crucial in the proof of the Cauchy stability of Minkowski space when treating nonlinear terms. 

The Schwarzschild solution describes a static, non-rotating black hole with mass $M$. The exterior region of  Schwarzschild black hole is described by $(t,r,\theta,\phi)=(t,r,\omega)\in\mathbb{R}\times(2M,\infty)\times S^2$ with the metric $g=-H\mathrm{d}t^2+H^{-1}\mathrm{d}r^2 +r^2(\mathrm{d}\theta^2+sin^2\theta\mathrm{d}\phi^2)$ where $H=1-2M/r$. At $r=3M$, there are null geodesics which orbit the black hole. Since null geodesics are the characteristics for the wave equation, one expects that energy can remain concentrated there for arbitrarily long periods of time, potentially violating any Morawetz estimate. Using $T=\partial_t$ and some modifications of the vector field $A=(1-3M/r)\partial_r$, one can repeat the argument from Minkowski space to show\cite{BlueSterbenz,BlueSoffer:Original,BlueSoffer:Errata,DafermosRodnianski:RedShift}
\begin{eqnarray*}
E_T(t)&=&E_T(0)=\int \left(|\partial_tu|^2+H^2|\partial_r u|^2+\frac{|\nabla_{S^2}u|}{r^2}\right) r^2\mathrm{d}r\mathrm{d}^2\omega, \\
E_T(0)&\geq& C\int_{0}^\infty\int \left(\frac{H^2|\partial_r u|^2}{r^2}+(1-3M/r)^2\left(\frac{|\partial_t u|^2}{r^2}+\frac{|\nabla_{S^2}u|^2}{r^3}\right)\right)r^2\mathrm{d}r\mathrm{d}^2\omega\mathrm{d}t . 
\end{eqnarray*}
($\mathrm{d}^2\omega$ denotes $\sin\theta\mathrm{d}\theta\mathrm{d}\phi$.) For null geodesics, the corresponding quantity $e_A$ measures the radial momentum pointing away from the orbiting null geodesics. This is increasing, and underlies the Morawetz estimate. The degeneracy at $r=3M$ of this estimate is sufficient to allow for the slow dispersion near the orbiting geodesics. The Schwarzschild spacetime is spherically symmetric so that, using the spherical Laplacian, there is a strengthened energy that is also conserved. 

This brief note ignores several important issues. For example, the geometrically induced measure on a hypersurface of constant $t$ is $H^{-1/2}r^2\mathrm{d}r\mathrm{d}^2\omega$, and, since $\lim_{r\rightarrow2M}g(T,T)=0$, the integrand in $E_T$ is degenerate with respect to this weight (Note that it has a weight of only $r^2\mathrm{d}r\mathrm{d}^2\omega$ without $H^{-1/2}$). Other ignored issues include the extended Schwarzschild spacetime, causal structure, and particularly the red-shift effect. These ideas can be used with a Morawetz estimate to construct a bounded, non-degenerate energy in all cases considered in this note and a remarkably broad set of other cases \cite{DafermosRodnianski:RedShift,DafermosRodnianski:KerrEnergyBounds}. 


\section{More complicated examples}
\label{S:Kerr}
Kerr's two parameter family of solutions to the Einstein equations
describe rotating black holes with mass $M$ and angular momentum $Ma$,
when $|a|<M$ (in natural units). For the Maxwell equation in the
rotating Kerr spacetime there are several crucial obstacles that must
be overcome in trying to adapt the method of the wave equation in the
Schwarzschild spacetime:
\begin{enumerate}
\item There is no globally time-like Killing vector, and, hence, no positive, conserved energy. 
\item The structure of the orbiting null geodesics is much more complicated. In particular, it fills an open set in spacetime. 
\item There are insufficiently many (classical) symmetries to strengthen the energy sufficiently to control the $L^\infty$ norm. 
\item The Maxwell equation is a system, instead of a scalar equation. 
\item There are bound states. 
\end{enumerate}
The first three problems are present in the study of the wave equation in the Kerr spacetime. The last two are present in the study of the Maxwell field in the Schwarzschild spacetime and in the study of the Einstein equation in asymptotically Schwarzschildian spacetimes. 

The first problem is quite severe since a positive, conserved energy is the basis for the study of most hyperbolic equations. This problem can be overcome by generating an energy with a vector field $T_\chi$ that is globally time like in the exterior and which fails to be Killing only in a fixed region away from the trapping, such as $r\in(5M,6M)$. Such a vector exists when $a$ is sufficiently small. Although $E_{T_\chi}$ is not conserved, its growth is controlled by the terms in the Morawetz estimate with an additional factor $|a|$. Thus, a boot-strap argument provides a uniform estimate in $t$, $E_T(t)\leq C E_T(0)$ once a Morawetz estimate can be proved. 

Since the orbiting null geodesics fill an open set in the exterior of the black hole, it is not possible to construct a vector field $A$ that points away from the orbiting null geodesics. However, in the tangent space, the set of orbiting null geodesics has codimension $3$, so it is relatively easy to construct a quantity that measures the distance in the tangent space from the orbiting null geodesics. Remarkably, this function can be constructed purely from conserved quantities of the null geodesic and functions of $r$. In accordance with the expectation from the wave particle analogy, corresponding to the null geodesic conserved quantities, there are symmetry operators $S_{\underline{a}}$ and separation constants for the wave equation. The operators and separation constants correspond to the constants of motion for the geodesic equation, and (with L. Andersson) we have constructed\cite{AnderssonBlue} a collection of vector fields $A^{\underline{a}\underline{b}}$ such that the fifth-order differential operator $S_{\underline{a}}A^{\underline{a}\underline{b}}S_{\underline{b}}$ corresponds to the distance in the tangent space from the orbiting geodesics. We replaced the energy-momentum tensor by a bilinear form $\mathrm{T}[u,v]=(1/4)(\mathrm{T}[u+v]-\mathrm{T}[u-v])$, define momenta and energies generated by collections of vector fields $A^{\underline{a}\underline{b}\beta}$ by contracting against $\mathrm{T}[S_{\underline{a}}u,S_{\underline{b}}u]$. This $A^{\underline{a}\underline{b}}$ then gives the desired Morawetz estimate. This requires working with $H^3$ type energies, because the additional differential operators have been applied to $u$. Similar vector fields were constructed using the separation constants, which has the advantage of remaining at the level of $H^1$ regularity. The method of separation of variables was used to prove energy bounds first\cite{DafermosRodnianski:KerrEnergyBounds} and a proof of this with a Morawetz estimate for all modes was developed independently afterwards\cite{TataruTohaneanu}. Our independent work using commuting operators was completed somewhat later. Unavoidably, there is some degeneracy in all approaches near the orbitting null geodesics. 

The lack of standard symmetry operators is compensated for by the existence of additional, conserved quantities for the null geodesics. As with the total angular momentum squared in $\mathbb{R}^{1+3}$, these are quadratic in $\dot\gamma$, but, unlike in $\mathbb{R}^{1+3}$, these do not admit a decomposition into Killing vectors. In addition to the quantities coming from the $2$ Killing vectors, there is an additional quantity that comes from an (irreducible) Killing $2$-tensor ($K_{\alpha\beta}=K_{(\alpha\beta)}$, $\nabla_{(\gamma}K_{\alpha\beta)}=0$.) The operator that cannot be decomposed into differentiation along Killing vectors is called a hidden symmetry. These were used in the construction of $A^{\underline{a}\underline{b}}$ and can also be used to obtain stronger Sobolev norms, although other methods exist\cite{DafermosRodnianski:KerrEnergyBounds}. 

Turning to the problems arising in the Maxwell and Einstein equations, it is somewhat surprising that the wave particle analogy fails to give a useful Morawetz estimate for these systems. Bound states exists, and these must violate any putative decay estimate. However, both systems can be decomposed into (complex) scalar components, $\phi_{-1}, \phi_0, \phi_1$ for the Maxwell field and $\psi_{-2},\psi_{-1},\psi_{0},\psi_{1},\psi_{2}$ for the (linearised or true) Einstein equations. 

In the Schwarzschild case, it has long been known that it is possible to derive a decoupled equation for the $\phi_0$ and $\psi_0$ components\cite{ReggeWheeler,PriceI,PriceII}. These equations are of the form $(\nabla^\alpha\nabla_\alpha+s^2V)(r^s\varphi_0)=0$, where $s=1$ ($s=2$) for the Maxwell (Einstein) equation, $V=2Mr^{-3}$ is a real potential, and $\varphi$ is $\phi$ or $\psi$ depending on the equation being reduced. For these wave equations with potential, one can apply the standard wave analysis and, fortunately, project out the bound states, and derive a Morawetz estimate for $\varphi_0$. This estimate is then sufficiently strong to yield decay estimates for the remaining components of the relevant field equation. (For the Maxwell equation in Schwarzschild, the bound states are spherically symmetric, correspond to the charge, and can be easily subtracted off from the initial data\cite{Blue:Maxwell}. For the Einstein equation in asymptotically Schwarzschildean spacetimes, the vanishing of the bound states is taken as an assumption\cite{Holzegel:AsymptoticallySchwarzschild}.) 

In going to the Maxwell equation in the Kerr spacetime (and similarly, for the linearised Einstein equations in the Kerr spacetime or, perhaps, the true Einstein equations in an asymptotically Kerr spacetime), one could attempt the same approach. It has long been known that the extreme components of each field $\varphi_{\pm s}$ satisfies equations famously derived by Teukolsky\cite{Teukolsky}. However, these are second-order equations, and, when put in the form $(\partial_\alpha L^{\alpha\beta}(s)\partial_\beta +W(s))(\tilde{\varphi}_{\pm s})=0$, the matrix of coefficients $L$ is not symmetric, which prevents most of the standard tools of hyperbolic PDE from being applied. The middle Maxwell component $\phi_0$ satisfies the equation\cite{FackerellIpser} $(\nabla^\alpha\nabla_\alpha +2M/p^3)(p\phi_0)=0$ with $p=r+ia\cos\theta$, and it has recently been shown that, for the linearised Einstein equation (in an appropriate gauge)\cite{AksteinerAndersson}, $(\nabla^\alpha\nabla_\alpha +8M/p^3)(p^2\phi_0)=0$. All of these equations can be separated. This suggests that there should be commuting operators, but, apart from the $\phi_0$ Maxwell component equation\cite{KalninsEtAl}, this does not appear to be in the literature. Unfortunately, it seems that the existence of a Killing tensor does not seem to be sufficient to guarantee the existence of a commuting operator, even when the Ricci curvature vanishes. 

Since the $\varphi_0$ wave equations have complex potentials, they lack a divergence-free energy-momentum tensor. Thus, when one construct a $T_{\chi}$ based energy, the new error terms that arise include terms of the form $\textrm{Im}(\bar{u}\partial_t u)$ with weights which do not vanish near the orbiting null geodesics. A model problem for this is the case of a wave equation with a complex potential on a static manifold with orbiting null geodesics occurring only at one radial value. This has recently been treated \cite{ABN:ModelProblem} by using pseudodifferential refinements of the Morawetz vector field, which allow for control over additional (fractional) derivatives in the Morawetz estimate, and which allow for control of $\textrm{Im}(\bar{u}\partial_t u)$ near the orbiting null geodesics. 

This suggests that it should be possible to treat the $\phi_0$ equation for the Maxwell equation similarly, although this will require a technically more cumbersome construction to construct a vector field from the separation constants and to then refine it using fractional powers of the separation constants. 

Since the proofs of the stability of Minkowski space have relied on energy-generation and energy-strengthening techniques, we would prefer to prove decay estimates for the Maxwell equation without resorting to separation of variables. To do this, with L. Anderson and J.P. Nicolas, we are currently considering approaches that avoid the decoupled $\phi_0$ equation. 

Other open problems inspired by this include: 
\begin{enumerate}
\item Classical analysis of the wave equation with a complex potential and trapping: Is it possible to estimate $\int \textrm{Im}(\bar{u}\partial_t u)r^{-k} \mathrm{d}x\mathrm{d}t$, for sufficiently large $k$, without using pseudodifferential techniques? 
\item The equation for $\phi_0$ can be rewritten as $\nabla_\alpha p^{-1}g^{\alpha\beta}\partial_\beta u=0$ with $u=p^2 \phi_0$. Since for $a$ small, $p$ seems very close to $r$, in some sense $p^{-1}g$ is almost a Lorentzian metric. Can such ``almost hyperbolic'' problems be treated? 
\item Can non-symmetric problems, like the Teukolsky equations, be treated using energy methods? 
\end{enumerate}

\section*{Acknowledgement}
I would like to thank my collaborators on this project, Lars Andersson and Jean-Philippe Nicolas. I would also like to thank the organisers of the ICMP and of the general relativity sessions. 

\bibliographystyle{ws-procs975x65}
\bibliography{blue}

\end{document}